\documentclass[a4paper,11pt,reqno]{article}
\usepackage[numbers,sort&compress]{natbib}
\usepackage{amsmath, amsfonts, amssymb, amsthm}
\input{epsf.sty}
\usepackage{mathrsfs}
\usepackage{amsmath}
\usepackage{amssymb}
\usepackage{graphicx}
\usepackage{subfigure}

\def \R {{\mathbb{R}}}
\numberwithin{equation}{section}

\begin{document}

\title{Stability of the train of $N$ solitary waves for the two-component Camassa-Holm shallow water system}

\author{Xingxing Liu\footnote{E-mail address:
liuxxmaths@cumt.edu.cn and tel: +86 15862186482}\\
Department of mathematics, China University of Mining and
Technology, \\Xuzhou, Jiangsu, 221116 China}
\date{}
\maketitle

\maketitle
\begin{abstract}
Considered herein is the integrable two-component Camassa-Holm shallow water system derived in the context of shallow water theory,
which admits blow-up solutions and the solitary waves interacting like solitons. Using modulation theory, and combining the almost monotonicity of a local version of energy
with the argument on the stability of a single solitary wave, we prove that the train of $N$
solitary waves, which are sufficiently decoupled, is orbitally stable in the energy space $H^1(\R)\times
L^2(\R)$.


\end{abstract}

\section{Introduction}
\newtheorem {remark1}{Remark}[section]
\newtheorem{theorem1}{Theorem}[section]
\newtheorem{proposition1}{Proposition}[section]
\newtheorem{lemma1}{Lemma}[section]
In this paper, we are concerned with the following two-component Camassa-Holm
shallow water system \cite{Chen,Constantin8,Ivanov,Popo}
\begin{equation}\label{2-CH}
\left\{\begin{array}{ll}m_t+2 u_x m +u m_x+\rho{\rho}_x=0,\ m=u-u_{xx}, & t>0,x\in \mathbb{R},\\
\rho_{t}+(u\rho)_x=0,  & t>0,x\in \mathbb{R},
\end{array}\right.
\end{equation}
where the variables $u(t,x)$, $\rho(t,x)$ describe the horizontal velocity of the fluid and the horizontal
deviation of the surface from equilibrium (or scalar density), respectively. The system (\ref{2-CH}) was originally introduced
by Chen et al. \cite{Chen} and Falqui in \cite{Falqui}. It is completely integrable \cite{Constantin8,Falqui,Ivanov1} as
it can be written as a compatibility condition of two linear systems (Lax pair). Compared with the other integrable multicomponent
Camassa-Holm-type systems, the system (\ref{2-CH}) has caught a large amount of attention, after Constantin and Ivanov
\cite{Constantin8} derived it in the context of shallow water regime. It is noticed that
the boundary assumptions $u\rightarrow 0$ and $\rho \rightarrow 1$ as $|x|\rightarrow \infty$, at any instant $t$, are required in their hydrodynamical derivation.

For $\rho\equiv0$, the system (\ref{2-CH}) becomes the classical
Camassa-Holm equation, which was first derived as an abstract bi-Hamiltonian partial differential equation by Fokas and Fuchssteiner \cite {FF}.
Then Camassa and Holm \cite{C-H} independently rediscovered it modeling shallow water waves with $u(t,x)$ representing the free surface
over a flat bottom. Moreover, it was found by Dai \cite{Dai} as a model for nonlinear waves in cylindrical hyperelastic rods where
$u(t,x)$ stands for the radial stretch relative to a pre-stressed state. In the past two decades, the reason for the Camassa-Holm equation as the master equation
in shallow water theory is that it gives a positive response to the question 'What mathematical models for shallow water waves could include both the phenomena of
soliton interaction and wave breaking?', which was proposed by Whitham \cite{Whitham}. The appearance of breaking waves as one of the remarkable properties of the
Camassa-Holm equation, however, can not be captured by the KdV and BBM equations \cite{C-L}. Plenty of impressive known results on wave breaking for the Camassa-Holm equation
have been obtained in \cite{C-H1,Constantin1,Constantin5,Constantin6,L-O,McKean,Molinet,Zhou}. Recently, we notice that Brandolese \cite{Brandolese1} unifies some of earlier results
by a more natural blow-up condition, that is, local-in-space blow-up criterion which means the condition on the initial data is purely local in space variable.

On the other hand, it was shown that the Camassa-Holm equation has solitary waves interacting like solitons \cite{C-H,C-H1}, which capture the essential
features of the extreme water waves \cite{Con,Con1,Con-J,Toland}.
Hence many papers addressed another fundamental qualitative property of solutions for the Camassa-Holm equation, which is the stability of solitary wave solutions.
As commented in \cite{C-S}, due to the fact that a small perturbation of a solitary wave can yield another one with a different speed and phase shift, we could only expect orbital stability
for solitary waves. Constantin and Strauss \cite{C-S} gave a very simple proof of the orbital stability of the peakons by using the conservation laws.
Then they \cite{C-S4} applied the general approach developed by \cite{G-S-S}, to cope with the stability of the smooth solitary waves. A series of works by
El Dika and Molinet \cite{E-M1,E-M2,E-M3} were devoted to the study of the stability of the train of $N$ solitary waves, multipeakons and multi antipeakon-peakons, respectively.
Moreover, Lenells \cite{Lenells} presented a variational proof of the stability of the periodic peakons.

For $\rho\neq0$, the system (\ref{2-CH}) has also attracted much attention
owing the fact that it has both solutions which blow up in finite time and solitary wave solutions interacting like solitons.
The Cauchy problem of the system (\ref{2-CH}) has been studied extensively. The local
well-posedness for the system (\ref{2-CH}) with initial data $(u_0,\rho_0)^t\in H^s\times H^{s-1},s\geq 2,$
by Kato's semigroup theory \cite{K1}, was established in \cite{Escher}. Then Gui and Liu \cite{Gui} improved the well-posedness
result with initial data in the Besov spaces (especially in $H^s\times H^{s-1}, s>\frac{3}{2}$). More interestingly,
singularities of the solutions for the system (\ref{2-CH}) can occur only in the form of wave breaking, while blow-up
solutions with a different class of certain initial profiles were shown in \cite{Constantin8,Escher,Guan1,Gui1,Gui,Z-L}.
Moreover, the system (\ref{2-CH})
has also global strong solutions \cite{Constantin8,Guan1,Gui1}. Here we recall the following global existence result needed in our developments.
\begin{proposition1}\label{Pro1}
(Global existence) Assume $\vec{u}^t_0(0,x)=(u_0(x),\eta_0(x))^t\in
H^s(\R)\times H^{s-1}(\R),s\geq2$, and $T>0$ be the maximal time of
the existence of the solution $\vec{u}^t(t,x)=\big(u(t,x),\eta(t,x)\big)^t\in
C\big([0,T);H^s(\R)\times H^{s-1}(\R)\big)\bigcap
C^1\big([0,T);H^{s-1}(\R)\times H^{s-2}(\R)\big)$ to the system
(\ref{2-CH}) with the initial profile $\vec{u}^t(0,x)$. If
$\eta_0(x)\neq -1$, then $T=+\infty$, $i.e.,$ the solution
$\vec{u}^t$ of the system (\ref{2-CH}) is global.
\end{proposition1}
On the other hand, Mustafa \cite{Mustafa} proved that the existence
of the smooth solitary waves for the system (\ref{2-CH}) with a single crest profile of maximum amplitude. For convenience,
we briefly collect the properties of the solitary waves from \cite{Mustafa,Z-L} in the proposition below.
\begin{proposition1}\label{Pro2}
(Existence of solitary waves) The system
(\ref{2-CH}) admits the smooth solitary wave solution
$\vec{\varphi}^t_c(t,x)\triangleq(\varphi_c(t,x), \xi_c(t,x))^t$ with the
speed $c>1$. Moreover, as $|x|\rightarrow \infty$, we have
\begin{eqnarray}\label{1.2}
\varphi(x)=O\big(exp(-\sqrt{\frac{c^2-1}{c^2}}|x|)\big),
\end{eqnarray}
and thus $\xi(x)=\frac{\varphi}{c-\varphi}$ also holds true.
\end{proposition1}
Recently, Zhang and Liu \cite{Z-L} obtained the orbital stability of a single solitary wave of the system (\ref{2-CH}) by following the general
spectral method developed by Grillakis et al. \cite{G-S-S}. Thus, an interesting problem is to investigate whether or not the train of
$N$ solitary wave solutions of the system (\ref{2-CH}) is orbitally stable as the scalar Camassa-Holm equation \cite{E-M1}. This is
the question we shall discuss in our present paper. For this purpose, we firstly rewrite the system (\ref{2-CH}) with $\rho\triangleq1+\eta$
($\eta\rightarrow 0$ as $|x|\rightarrow\infty$) as follows
\begin{equation}\label{2-CH2}
\left\{\begin{array}{ll}u_t-u_{txx}+3uu_x=2u_xu_{xx}+uu_{xxx}-(1+\eta)\eta_x, & t>0,x\in \mathbb{R},\\
\eta_{t}+(u(1+\eta))_x=0,  & t>0,x\in \mathbb{R},
\end{array}\right.
\end{equation}
or equivalently,
\begin{equation}\label{2-CH3}
\left\{\begin{array}{ll}(1-\partial^2_x)u_t=-\frac{1}{2}(1-\partial^2_x)\partial_xu^2
-\partial_x(u^2+\frac{1}{2}u^2_x+\eta+\frac{1}{2}\eta^2)
, & t>0,x\in \mathbb{R},\\
\eta_{t}=-\partial_x(u+u\eta),  & t>0,x\in \mathbb{R}.
\end{array}\right.
\end{equation}
Secondly, we define the space $X=H^1(\R)\times L^2(\R)$ with inner product $(\cdot,\cdot)$ and its norm $\|\cdot\|_X$, and thus the
dual $X^*=H^{-1}(\R)\times L^2(\R)$. Denote by $\langle\cdot,\cdot\rangle$ the pairing between $X$ and $X^*$, and the space $X^{**}$ is identified with $X$ in the natural way.
We have the isomorphism $I:X\rightarrow X^*$ defined by
\begin{equation*}I=\left(
\begin{array}{cc}
1-\partial_x^2 & 0\\
0 & 1\\
\end{array}
\right).\end{equation*} Thus, for $\vec{u}^t,\vec{v}^t\in X$, we obtain $\langle I\vec{u}^t,\vec{v}^t\rangle=(\vec{u}^t,\vec{v}^t)$. Moreover, we also need the functionals on $X$, which
are two useful conservation laws \cite{Constantin8,Z-L}
\begin{equation}\label{1.5}
E(\vec{u}^t)\triangleq \frac{1}{2}\int_{\R}(u^2+u^2_x+\eta^2) dx  \quad
\mbox{and}\quad F(\vec{u}^t)\triangleq\frac{1}{2}\int_{\R}
(u^3+uu_x^2+2u\eta+u\eta^2)dx.
\end{equation}
Now we are in the position to state our main result.
\begin{theorem1}\label{Th1}
Given $N$ velocities $c_1,...,c_N$ such that $1<c_1<c_2<...<c_N$.
Denote the $N$ solitary waves as $\sum\limits_{i=1}\limits^N
\vec{\varphi}^t_{c_i}(t,x)=\sum\limits_{i=1}\limits^N
(\varphi_{c_i}(x-c_it),\xi_{c_i}(x-c_it))^t$ by Proposition \ref{Pro2}. There
exist $\gamma_0,A_0,L_0,\varepsilon_0>0$ such that for any initial
date $\vec{u}^t_0=\vec{u}^t(0,x)=(u_0,\eta_0)^t\in H^s(\R)\times
H^{s-1}(\R),$ with $s\geq 2$, $\eta_0(x)\neq -1$ and satisfying
\begin{eqnarray}\label{1.6}
\big\|\vec{u}^t_0-\sum_{i=1}^N\vec{\varphi}^t_{c_i}(\cdot-x_i)\big\|_{X}
\leq\varepsilon,
\end{eqnarray}
for some $0<\varepsilon<\varepsilon_0$, and $x_i-x_{i-1}\geq L$ with
$L>L_0$ for $i=2,...,N$. Then, for the global strong solution
$\vec{u}^t=(u,\eta)^t$ of the system (\ref{2-CH2}) or (\ref{2-CH3})  with
$\vec{u}^t_0$ guaranteed by Proposition \ref{Pro1}, there exist
$x_1(t),...,x_N(t)$ such that
\begin{eqnarray*}
\quad \quad
\sup_{t\in[0,\infty)}\big\|\vec{u}^t(t,\cdot)-\sum_{i=1}^N\vec{\varphi}^t_{c_i}(\cdot-x_i(t))\big\|_{X}
\leq A_0(\sqrt\varepsilon+e^{-\gamma_0L}).
\end{eqnarray*}
\end{theorem1}
As commented by El Dika and Molinet \cite{E-M1,E-M2}, the strategy initiated in \cite{M-M-T} for a scalar equation
indicates that there are principally two required ingredients to prove the stability of the sum of $N$ solitary waves. One of them is a
dynamical proof of the stability of the single solitary wave, and the other is a property of almost monotonicity, which says
for a solution close to $\vec{\varphi}^t_{c}$, the part of the energy traveling at the right of $\vec{\varphi}^t_{c}(\cdot-ct)$ is
almost time decreasing. Our approach to prove Theorem \ref{Th1} is try to follow this method. However, we consider here a coupled system
with two component, instead of a scalar equation. Hence, the same argument as in \cite{E-M1} or \cite{M-M-T} for a single equation is not
directly applicable here. More precisely, we need to overcome a difficulty encountered by the coupled system (\ref{2-CH}) in comparison with
the Camassa-Holm equation, which is the mutual effect between the two component $u(t,x)$ and $\eta(t,x)$.
To solve this problem, here we require more elaborate analysis on the decomposition of the solution by using modulation theory, and a local coerciveness
inequality related to a Hessian operator $H_c$ of $cE-F$ around the solitary wave $\vec{\varphi}^t_{c}$. Moreover, we know that
the method of the proof of the stability relies heavily on a property of almost monotonicity. Therefore, the key issue to prove Theorem \ref{Th1} is to
estimate precisely these coupled terms, which appear in the energy at the right of the $(j-1)$-th bump of the solution. To this end, we first
apply H\"{o}lder inequality to prove $(p\ast \eta)^2\leq p\ast \eta^2$ $(p\triangleq \frac{1}{2}e^{-|x|})$, and then deduce the desired result by means of
the Minkowski inequality (see Lemma \ref{Lem2} below). Actually, to the best of our knowledge, our theorem is the first result on the stability of the sum of $N$ solitary wave solutions
for the coupled shallow water system. Hence we expect there are more applications of this method to handle the stability of $N$ solitary waves
for the other two-component system, such as a generalized two-component Camassa-Holm system \cite{C-L-Q}, two-component Dullin-Gottwald-Holm system \cite{Liu1}
and so on.

The remainder of the paper is dedicated to the proof of Theorem \ref{Th1}.
In Section 2, we present some useful lemmas which will be used in
the sequel. First, we control the distance between the different bumps of the solution by using a modulation argument.
Then, we prove a almost monotonicity property and local coercivity of the solitary wave. In Section 3, we complete the proof of Theorem
\ref{Th1} by three steps.

$Notation.$ As above and henceforth, we denote by $\ast$ the convolution.
Since our discussion is all on the line
$\R$, for simplicity, we omit $\R$ in our notations of function
spaces. All the transpose of a row vector $\vec{f}=(f_1,f_2)$ is
presented as $\vec{f}^t =\left(
\begin{array}{c}
f_1 \\
f_2\\
\end{array}
\right)$.

\section{Preliminaries}\label{Sec2}
\newtheorem {remark2}{Remark}[section]
\newtheorem{theorem2}{Theorem}[section]
\newtheorem{lemma2}{Lemma}[section]
In this section, we shall establish some useful lemmas which are
crucial to pursue our goal. We break them into the following three
subsections.
\subsection{Modulation}
In this subsection, we show that we can decompose the solution
$\vec{u}^t(t)=(u,\eta)^t$ as the sum of $N$ modulated solitary waves
and a vector function $\vec{v}^t(t)\triangleq(v,\zeta)^t$, which remains small in
the space $X$:
\begin{eqnarray*}
\vec{u}^t(t,x)=\sum_{i=1}^N\big(\varphi_{c_i}(x-x_i(t)),
\xi_{c_i}(x-x_i(t))\big)^t+\vec{v}^t(t,x),
\end{eqnarray*}
with $\vec{v}^t(t,x)$ is orthogonal to
$\big((1-\partial^2_x)\partial_x\varphi_{c_i}(x-x_i(t)),
\partial_x\xi_{c_i}(x-x_i(t))\big)^t$ in $L^2$, for $i=1,...,N.$
Moreover, we prove that the different bumps of $\vec{u}^t(t)$
that are individually close to a solitary wave get away from each
other as time is increasing.

Let $1<c_1<c_2...<c_N$,
$\sigma_0=\frac{1}{4}\min(c_1,c_2-c_1,...,c_N-c_{N-1})$. For
$\alpha, \ L >0$, we define the neighborhood of size $\alpha$ of the
superpositon of $N$ solitary waves of speed $c_i$, located at a
distance larger than $L$,
\begin{eqnarray*}
U(\alpha,L)=\{\vec{u}^t\in X;\
\inf_{x_i-x_{i-1}>L}\big\|\vec{u}^t-\sum_{i=1}^N\vec{\varphi}^t_{c_i}
(\cdot-x_i)\big\|_{X}<\alpha \}.
\end{eqnarray*}

\begin{lemma2}\label{Lem1}
Suppose $\vec{u}^t_0$ satisfy the assumptions (\ref{1.6}) given in Theorem \ref{Th1}.
There exist $\alpha_0,L_0$ such that for all $0<\alpha<\alpha_0$ and
$L>L_0$, if the corresponding solution $\vec{u}^t(t)\in
U(\alpha,\frac{L}{2})$ on $[0,t_0]$ for some $0<t_0\leq+\infty$,
then there exist unique $C^1$-functions $x_i:\ [0,t_0]\rightarrow
\R,i=1,...,N,$ such that
\begin{eqnarray*}
\vec{v}^t(t,x)=\big(v(t,x),\zeta(t,x)\big)^t=\vec{u}^t(t,x)-\sum_{i=1}^N\vec{R}^t_i(t,x),
\end{eqnarray*}
where
$\vec{R}^t_i(t,x)=\big(R_i(t,x),S_i(t,x)\big)^t=\big(\varphi_{c_i}(x-x_i(t)),
\xi_{c_i}(x-x_i(t))\big)^t $, satisfies the following orthogonality
conditions
\begin{eqnarray}\label{2.1}
\int_\R v(t)(1-\partial_x^2)\partial_xR_i(t)dx+\int_\R
\zeta(t)\partial_xS_i(t)dx=0,\quad \quad i=1,...,N.
\end{eqnarray}
Moreover, the following statements hold true:
\begin{eqnarray}\label{2.2}
\|\vec{v}^t(t)\|_{X}=\big\|\vec{u}^t(t)-\sum_{i=1}^N\vec{R}^t_i(t)\big\|_{X}=O(\alpha),
\end{eqnarray}
\begin{eqnarray}\label{2.3}
\sup\limits_{t\in[0,t_0]}|\dot{x}_i(t)-c_i|\leq
O(\alpha)+O(e^{-\sigma_0L}), \quad \quad i=1,...,N,
\end{eqnarray}
\mbox{and}
\begin{eqnarray}\label{2.4}
x_i(t)-x_{i-1}(t)\geq\frac{3L}{4}+2\sigma_0t, \quad \quad i=2,...,N.
\end{eqnarray}
\end{lemma2}
\begin{proof}
For $Z=(z_1,z_2,...,z_N)\in {\R}^N$, such that
$z_i-z_{i-1}>\frac{L}{2}$, we set
$\vec{R}^t_Z(\cdot)=(R_Z(\cdot),S_Z(\cdot))^t
=\big(\sum\limits_{i=1}\limits^N\varphi_{c_i}(\cdot-z_i),\sum\limits_{i=1}\limits^N\xi_{c_i}(\cdot-z_i)\big)^t
,$ and denote $B_{H^1}(R_Z,\alpha),B_{L^2}(S_Z,\alpha)$ as the ball
in $H^1,L^2$ of center $R_Z,S_Z$ with radius $\alpha$, respectively. For $0<\alpha<\alpha_0$, we define the following mapping
\begin{eqnarray*}
   Y:(-\alpha,\alpha)^N\times
B_{H^1}(R_Z,\alpha)\times B_{L^2}(S_Z,\alpha)
&\rightarrow& \R^N,\\
(y_1,...,y_N,u,\eta)&\mapsto&
\big(Y^1(y_1,...,y_N,u,\eta),...,Y^N(y_1,...,y_N,u,\eta)\big),
\end{eqnarray*}
with \begin{eqnarray}\label{2.5}
 Y^i(y_1,...,y_N,u,\eta)
&=&\int_{\R}\big(u-\sum\limits_{i=1}\limits^N
\varphi_{c_i}(\cdot-z_i-y_i)
\big)(1-\partial_x^2)\partial_x\varphi_{c_i}(\cdot-z_i-y_i)dx\nonumber\\
&&+\int_{\R}\big(\eta-\sum\limits_{i=1}\limits^N
\xi_{c_i}(\cdot-z_i-y_i) \big)\partial_x\xi_{c_i}(\cdot-z_i-y_i)dx.
\end{eqnarray}
In the following, we verify that the function $Y$ satisfies the properties:

(i)\ $Y(0,...,0,R_Z,S_Z)=(0,...,0).$

(ii)\ For $i,j=1,...,N$, by the dominated convergence theorem and
the smoothness of $\vec{\varphi}^t_c$, the partial derivatives
$\frac{\partial Y^i}{\partial y_j},\frac{\partial Y^i}{\partial u}$
and $\frac{\partial Y^i}{\partial \eta}$ are continuous. Indeed, for
$i=1,...,N$,
\begin{eqnarray*}
\frac{\partial Y^i}{\partial
y_i}(y_1,...,y_N,u,\eta)&=&\int_{\R}\big(u_x-\sum\limits_{j=1,j\neq
i}\limits^N\partial_x\varphi_{c_j}(\cdot-z_j-y_j)\big)(1-\partial^2_x)\partial_x\varphi_{c_i}(\cdot-z_j-y_j)dx
\\
&&+\int_{\R}\big(\eta_x-\sum\limits_{j=1,j\neq
i}\limits^N\partial_x\xi_{c_j}(\cdot-z_j-y_j)\big)\partial_x\xi_{c_i}(\cdot-z_j-y_j)dx,
\end{eqnarray*}
\begin{eqnarray*}
\frac{\partial Y^i}{\partial
y_j}(y_1,...,y_N,u,\eta)&=&\int_{\R}\partial_x\varphi_{c_j}(\cdot-z_j-y_j)(1-\partial^2_x)\partial_x\varphi_{c_i}(\cdot-z_j-y_j)dx
\\
&&+\int_{\R}\partial_x\xi_{c_j}(\cdot-z_j-y_j)\partial_x\xi_{c_i}(\cdot-z_j-y_j)dx,\
\mbox{for}\ \forall j\neq i,
\end{eqnarray*}
\begin{eqnarray*}
\frac{\partial Y^i}{\partial
u}(y_1,...,y_N,u,\eta)=\int_{\R}(1-\partial^2_x)\partial_x\varphi_{c_i}(\cdot-z_j-y_j)dx,
\end{eqnarray*}
and
\begin{eqnarray*}
\frac{\partial Y^i}{\partial
\eta}(y_1,...,y_N,u,\eta)=\int_{\R}\partial_x\xi_{c_i}(\cdot-z_j-y_j)dx.
\end{eqnarray*}
Thus, the function $Y$ is of class $C^1$.

(iii)\ The determinant of the matrix
$Y'_{(y_1,...,y_N)}(0,...,0,R_Z,S_Z)$ is not equal to zero. In fact,
from the above (ii), we obtain
\begin{eqnarray*}
&&\frac{\partial Y^i}{\partial
y_i}(0,...,0,R_Z,S_Z)\\
&=&\|\partial_x\varphi_{c_i}\|^2_{H^1}+\|\partial_x\xi_{c_i}\|^2_{L^2}
\geq
\frac{1}{2}\min\limits_{i=1,...,N}\big(\|\partial_x\varphi_{c_i}\|^2_{H^1}+\|\partial_x\xi_{c_i}\|^2_{L^2}\big),
\end{eqnarray*}
and for $j\neq i$, by the exponential decay (\ref{1.2}) of $\vec{\varphi}_c^t$,
and $z_i-z_{i-1}>\frac{L}{2}$, for $L_0$ large enough,
\begin{eqnarray*}
&&\frac{\partial Y^i}{\partial
y_j}(0,...,0,R_Z,S_Z)\\
&=&\big(\partial_x\varphi_{c_j}(x-z_j),\partial_x\varphi_{c_i}(x-z_i)\big)_{H^1}
+\big(\partial_x\xi_{c_j}(x-z_j),\partial_x\xi_{c_i}(x-z_i)\big)_{L^2}\leq
O(e^{-\sigma_0L}).
\end{eqnarray*}
Thus, we deduce that for $L_0$ large enough, $\frac{\partial
Y^i}{\partial y_i}(0,...,0,R_Z,S_Z)\gg \frac{\partial Y^i}{\partial
y_j} (0,...,0,R_Z,S_Z)$ whenever $j\neq i$. So,
$Y'_{(y_1,...,y_N)}(0,...,0,R_Z,S_Z)$ is a diagonally dominant
matrix.

Therefore, from the implicit function theorem, we find that there
exists $\alpha_0>0$ and the uniquely determined $C^1$-functions
$(y_1,...,y_N)$ from $B_{H^1}(R_Z,\alpha_0)\times
B_{L^2}(S_Z,\alpha_0)$ to a neighborhood of $(0,...,0)$, such that
$$Y(y_1,...,y_N,u,\eta)=(0,...,0),\quad \mbox{for\ all}\ \vec{u}^t\in B_{H^1}(R_Z,\alpha_0)\times
B_{L^2}(S_Z,\alpha_0).$$ In particular, if $\vec{u}^t\in
B_{H^1}(R_Z,\alpha)\times B_{L^2}(S_Z,\alpha)$ with $0<\alpha\leq
\alpha_0$, there exists a constant $C_0>0$ such that
\begin{eqnarray}\label{2.6}
\sum\limits_{i=1}\limits^N\big|y_i(\vec{u}^t)\big|\leq C_0\alpha.
\end{eqnarray}
Note that $\alpha_0$ and $C_0$ depend only on the velocity $c_1$ and
$L_0$ and not on the point $(z_1,...,z_N)$. For $\vec{u}^t\in
B_{H^1}(R_Z,\alpha_0)\times B_{L^2}(S_Z,\alpha_0)$, we set
$x_i(\vec{u}^t)= z_i+y_i(\vec{u}^t).$ If we take
$\alpha_0\leq\frac{L_0}{8C_0}$, then $(x_1,...,x_N)$ are
$C^1$-functions on $\vec{u}^t\in B_{H^1}(R_Z,\alpha)\times
B_{L^2}(S_Z,\alpha)$, satisfying
\begin{eqnarray}\label{2.7}
x_i(\vec{u}^t)-x_{i-1}(\vec{u}^t)>\frac{L}{2}-2C_0\alpha>\frac{L}{4}.
\end{eqnarray}
Then by a modulation argument and the construction (\ref{2.5}) of the
functions $Y^i$, we can define $N$ $C^1$-functions $t\mapsto
x_i(t)=x_i(\vec{u}^t(t))$ satisfying the orthogonality conditions
(\ref{2.1}) for $i=1,...,N$. Furthermore, from $\vec{u}^t(t)\in
U(\alpha,\frac{L}{2})$, (\ref{2.6}) and the triangular inequality, there
exists $C_0>0$, such that for all $t\in[0,t_0]$
\begin{eqnarray}\label{2.8}
\|\vec{v}^t(t)\|_{X}\leq C_0\alpha,\end{eqnarray}
hence, (\ref{2.2}) holds true.

Differentiating the orthogonality conditions (\ref{2.1}) with
respect to time $t$, we obtain
\begin{eqnarray*}
\int_\R v_t(1-\partial_x^2)\partial_xR_idx+\int_\R
\zeta_t\partial_xS_idx=\dot{x}_i\big(\int_\R
v(1-\partial_x^2)\partial^2_xR_idx+\int_\R\zeta\partial^2_xS_idx\big),
\end{eqnarray*}
and thus, we have
\begin{eqnarray}\label{2.9}
\big|\int_\R v_t(1-\partial_x^2)\partial_xR_idx+\int_\R
\zeta_t\partial_xS_idx\big|\leq|\dot{x}_i|O(\|\vec{v}^t\|_{X})\leq|\dot{x}_i
-c_i|O(\|\vec{v}^t\|_{X})+O(\|\vec{v}^t\|_{X}).
\end{eqnarray}
Substituting $u=v+\sum\limits_{i=1}\limits^NR_i$ and
$\eta=\zeta+\sum\limits_{i=1}\limits^NS_i$ into the system (\ref{2-CH3}), it
follows
\begin{eqnarray}\label{2.10}
(1-\partial_x^2)v_t+\sum\limits_{i=1}\limits^N(1-\partial_x^2)\partial_tR_i
&=&-\frac{1}{2}(1-\partial_x^2)\partial_x\big(
(v+\sum_{i=1}^NR_i)^2\big)-\partial_x\big((v+\sum_{i=1}^NR_i)^2\\
&&+\frac{1}{2}(v_x+\sum_{i=1}^N\partial_xR_i)^2+(\eta+\sum_{i=1}^NS_i)+
\frac{1}{2}(\eta+\sum_{i=1}^NS_i)^2 \big),\nonumber
\end{eqnarray}
and
\begin{eqnarray}\label{2.11}
\zeta_t+\sum\limits_{i=1}\limits^N\partial_tS_i
=-\partial_x\big((v+\sum_{i=1}^NR_i)
+(v+\sum_{i=1}^NR_i)(\eta+\sum_{i=1}^NS_i)\big).
\end{eqnarray}
From the definition of $R_i$ and $S_i$, we obtain
\begin{eqnarray}\label{2.12}
(1-\partial_x^2)\partial_tR_i+(\dot{x}_i-c_i)(1-\partial_x^2)\partial_xR_i
+3R_i\partial_xR_i=2\partial_xR_i\partial^2_xR_i+R_i\partial^3_xR_i-(1+S_i)\partial_xS_i,
\end{eqnarray}
and
\begin{eqnarray}\label{2.13}
\partial_tS_i+(\dot{x}_i-c_i)\partial_xS_i+\partial_xR_i+\partial_x(R_iS_i)=0.
\end{eqnarray}
Combining (\ref{2.10})-(\ref{2.13}), we infer that $v(t)$ satisfies on $[0,t_0]$
\begin{eqnarray}\label{2.14}
&&(1-\partial_x^2)v_t-\sum\limits_{i=1}\limits^N(\dot{x}_i-c_i)(1-\partial_x^2)\partial_xR_i\\
&=&-\frac{1}{2}(1-\partial_x^2)\partial_x\big(
(v+\sum_{i=1}^NR_i)^2-\sum_{i=1}^NR_i^2\big)
-\partial_x\big((v+\sum_{i=1}^NR_i)^2-\sum_{i=1}^NR_i^2+\frac{1}{2}(v_x+\sum_{i=1}^N\partial_xR_i)^2\nonumber\\
&&-\frac{1}{2}\sum_{i=1}^N(\partial_xR_i)^2 \big)-\partial_x\big(
\eta+ \frac{1}{2}(\eta+\sum_{i=1}^NS_i)^2
-\frac{1}{2}\sum_{i=1}^NS_i^2\big),\nonumber
\end{eqnarray}
and for $\zeta(t)$
\begin{eqnarray}\label{2.15}
\zeta_t-\sum\limits_{i=1}\limits^N(\dot{x}_i-c_i)\partial_xS_i
=-\partial_x\big(v
+(v+\sum_{i=1}^NR_i)(\eta+\sum_{i=1}^NS_i)-\sum_{i=1}^NR_iS_i\big).
\end{eqnarray}
Taking the $L^2$-scalar product (\ref{2.14}) with $\partial_xR_j$ and
(\ref{2.15}) with $\partial_xS_j$, integrating by parts, using the
exponential decay of $\vec{R}^t_i$ and its derivatives, by
(\ref{2.8})-(\ref{2.9}) and (\ref{2.7}), then plugging the two obtained results, we
get
\begin{eqnarray*}
|\dot{x}_j-c_j|\big(\|\partial_xR_j\|^2_{H^1}+\|\partial_xS_j\|^2_{L^2}+O(
\alpha)\big)\leq O(\alpha)+O(e^{-\sigma_0L}).
\end{eqnarray*}
Taking $\alpha_0$ small enough and $L_0$ large enough depending only on $\{c_i\}^N_{i=1}$, we obtain (\ref{2.3}).

To prove (\ref{2.4}), for $\alpha_0$ sufficient small and $L_0$ large
enough, we have $|\dot{x_i}-c_i|\leq\frac{c_i-c_{i-1}}{4}$. Thus for
all $0<\alpha<\alpha_0$ and $L\geq L_0>4C_0\varepsilon$, by the mean value theorem, (\ref{1.6}), (\ref{2.6}) and (\ref{2.3}),
there exist $\xi\in [0,t]$ such that
\begin{eqnarray*}
x_i(t)-x_{i-1}(t) &=&x_i(t)-x_i(0)+x_{i}(0)-x_{i-1}(0)
+x_{i}(0)-x_{i-1}(t)\nonumber\\
&=&x_{i}(0)-x_{i-1}(0)+(\dot{x}_{i}(\xi)-\dot{x}_{i-1}(\xi))t\nonumber\\
&>&L-C_0\varepsilon+\frac{(c_i-c_{i-1})t}{2}
\geq\frac{3L}{4}+2\sigma_0t,\nonumber \quad \forall\ t\in[0,t_0].
\end{eqnarray*}
This completes the proof of Lemma \ref{Lem1}.
\end{proof}
\subsection{Monotonicity property}
\par
This subsection is devoted to proving the principal tool of our proof,
which is the almost monotonicity of functionals that are very close
to the energy at the right of the $(i-1)$th bump of $\vec{u}^t$,
$i=2,...,N$. Firstly, we define $\Psi$ to be a $C^\infty$ function
such that
\begin{equation*}
\Psi(x)=\left\{\begin{array}{ll}e^{-|x|}, & x<-1,\\
1-e^{-|x|},& x>1,
\end{array}\right.
\ \mbox{and} \
\left\{\begin{array}{ll}0<\Psi\leq1,\Psi'>0,\quad   &x\in \mathbb{R},\\
 |\Psi'''|\leq10\Psi', \quad  &x\in [-1,1].\\
\end{array}\right.\end{equation*}
Consider $\Psi_K=\Psi(\frac{\cdot}{K})$, where the constant $K>0$
will be chosen later. Then, we introduce for $j=2,...,N,$
\begin{eqnarray*}
I_{j,K}(t)=I_{j,K}(t,\vec{u}^t(t))=\frac{1}{2}\int_\R
(u^2(t)+u^2_x(t)+\eta^2(t))\Psi_{j,K}(t,x)dx,
\end{eqnarray*}
where $\Psi_{j,K}(t,x)=\Psi_K(x-y_j(t))$ with
$y_j(t)\triangleq\frac{x_{j-1}(t)+x_j(t)}{2}$ for $j=2,...,N$. And for
$i=1,...,N,$ we define the following localized version of the
conservation laws (\ref{1.5}) of $E$ and $F$ as
\begin{eqnarray}\label{2.16}
E^t_i(\vec{u}^t)=\frac{1}{2}\int_\R \Phi_i(t)(u^2+u^2_x+\eta^2)dx\
\mbox{and} \ F^t_i(\vec{u}^t)=\frac{1}{2}\int_\R
\Phi_i(t)(u^3+uu^2_x+2u\eta+u\eta^2)dx,
\end{eqnarray}
where the weight functions $\Phi_i=\Phi_i(t,x)$ are given by
\begin{eqnarray*}
\Phi_1=1-\Psi_{2,K}, \ \Phi_N=\Psi_{N,K}\ \  \mbox{and}\ \
 \Phi_i=\Psi_{i,K}-\Psi_{{i+1},K},\ i=2,...,N-1.
\end{eqnarray*}
Obviously, $\sum\limits_{i=1}^N\Phi_{i,K}\equiv1$. Taking $L>0$ and
$\frac{L}{K}>0$ large enough, we can deduce that
\begin{eqnarray}\label{2.17}
\big|1-\Phi_{i,K}\big|\leq4e^{-\frac{L}{4k}}, \quad \mbox{for}\
x\in[x_i-\frac{L}{4},x_i+\frac{L}{4}],
\end{eqnarray}
and
\begin{eqnarray}\label{2.18}\quad\quad\quad
\big|\Phi_{i,K}\big|\leq4e^{-\frac{L}{4k}}, \quad \mbox{for}\
x\in[x_j-\frac{L}{4},x_j+\frac{L}{4}], \ \mbox{whenever}\ j\neq i.
\end{eqnarray}

\begin{lemma2}\label{Lem2}
Let $\vec{u}^t(t,x)$ be the global strong solution of the system
(\ref{2-CH}) such that $\vec{u}^t(t)\in U(\alpha,\frac{L}{2})$ on
$[0,+\infty)$, where $x_i(t)$ are defined in Lemma \ref{Lem1}. There exist
$\alpha_0>0$ and $L_0>0$ only depending on $\sigma_0$, such that if
 $0<\alpha<\alpha_0$ and $L\geq L_0$ then for any $5\leq K=O(\sqrt{L})$,
\begin{eqnarray}\label{2.19}
I_{j,K}(t)-I_{j,K}(0)\leq O(e^{-\sigma_0L}), \quad t\in[0,+\infty).
\end{eqnarray}
\end{lemma2}
\begin{proof}
By the conversation law $F(\vec{u}^t)$, we see that the system
(\ref{2-CH}) can be written in Hamiltonian form as
\begin{eqnarray}\label{2.20}
\frac{\partial \vec{u}^t}{\partial t}=JF'(\vec{u}^t),
\end{eqnarray}
where $J$ is a closed skew symmetric operator defined by
\begin{equation*}J=\left(
\begin{array}{cc}
-\partial_x(1-\partial_x^2)^{-1} & 0\\
0 & -\partial_x\\
\end{array}
\right),\end{equation*} and $F'(\vec{u}^t)$ is the Fr\'{e}chet
derivatives of $F$ in $X$ at $\vec{u}^t$, which can be calculated as
\begin{eqnarray}\label{2.21}
F'_u=\frac{3}{2}u^2-\frac{1}{2}u^2_x-uu_{xx}+\eta+\frac{1}{2}\eta^2\quad
\mbox{and}\quad F'_\eta=u+u\eta.
\end{eqnarray}

Differentiating $I_j(t)$ with respect to time $t$, we have
\begin{eqnarray}\label{2.22}
\frac{d}{dt}I_j(t)&=&-\frac{\dot{y}_j(t)}{2}\int_{\R}(u^2+u^2_x+\eta^2)\Psi'_{j,k}dx+
\int_{\R}(uu_t+u_xu_{xt}+\eta\eta_t)\Psi_{j,k}dx\nonumber\\
&=&-\frac{\dot{y}_j(t)}{2}\int_{\R}(u^2+u^2_x+\eta^2)\Psi'_{j,k}dx+J(t).
\end{eqnarray}
Integrating by parts, and by (\ref{2.20})-(\ref{2.21}), we deduce
\begin{eqnarray}\label{2.23}
J(t)&=&\int_{\R}uu_t\Psi_{j,k}dx-\int_{\R}uu_{txx}\Psi_{j,k}dx-\int_{\R}uu_{tx}\Psi'_{j,k}dx
+\int_{\R}\eta\eta_t\Psi_{j,k}dx\nonumber\\
&=&\int_{\R}u\Psi_{j,k}(1-\partial^2_x)u_tdx-\int_{\R}uu_{tx}\Psi'_{j,k}dx
+\int_{\R}\eta\eta_t\Psi_{j,k}dx\nonumber\\
&=&-\int_{\R}u\Psi_{j,k}\partial_xF'_udx+\int_{\R}u\Psi'_{j,k}(1-\partial^2_x)^{-1}\partial^2_xF'_udx
-\int_{\R}\eta\Psi_{j,k}\partial_xF'_\eta dx\nonumber\\
&=&\int_{\R}u_x\Psi_{j,k}F'_udx+\int_{\R}u\Psi'_{j,k}(1-\partial^2_x)^{-1}F'_udx
+\int_{\R}\eta_x\Psi_{j,k}F'_\eta dx+\int_{\R}\eta\Psi'_{j,k}F'_\eta
dx\nonumber\\
&\triangleq&J_1+J_2+J_3+J_4.
\end{eqnarray}
Integrating by parts and by (\ref{2.21}), we can directly compute
$J_i,1\leq j \leq 4$ as follows
\begin{eqnarray*}
J_1&=&\int_{\R}u_x\Psi_{j,k}(\frac{3}{2}u^2-\frac{1}{2}u^2_x-uu_{xx}+\eta+\frac{1}{2}\eta^2)dx\\
&=&-\frac{1}{2}\int_{\R}\Psi'_{j,k}u^3dx+\frac{1}{2}\int_{\R}\Psi'_{j,k}uu^2_xdx+\int_{\R}u_x\Psi_{j,k}\eta
dx+\frac{1}{2}\int_{\R}u_x\Psi_{j,k}\eta^2dx,
\end{eqnarray*}
\begin{eqnarray*}
J_2&=&\int_{\R}u\Psi'_{j,k}(1-\partial^2_x)^{-1}(\frac{3}{2}u^2-\frac{1}{2}u^2_x-uu_{xx}+\eta+\frac{1}{2}\eta^2)dx\\
&=&\int_{\R}u\Psi'_{j,k}(1-\partial^2_x)^{-1}(\frac{3}{2}u^2+\frac{1}{2}u^2_x-\frac{1}{2}(u^2)_{xx}+\eta
+\frac{1}{2}\eta^2)dx\\
&=&\frac{1}{2}\int_{\R}\Psi'_{j,k}u^3dx+\int_{\R}u\Psi'_{j,k}(1-\partial^2_x)^{-1}(u^2+\frac{1}{2}u^2_x+\eta
+\frac{1}{2}\eta^2)dx,
\end{eqnarray*}
\begin{eqnarray*}
J_3&=&\int_{\R}\eta_x\Psi_{j,k}(u+u\eta)
dx\\
&=&-\int_{\R}u_x\Psi_{j,k} \eta dx-\int_{\R}u\eta
\Psi'_{j,k}dx-\frac{1}{2}\int_{\R}u_x\Psi_{j,k} \eta^2
dx-\frac{1}{2}\int_{\R} u\Psi'_{j,k} \eta^2 dx,
\end{eqnarray*}
and
\begin{eqnarray*}
J_4&=&\int_{\R}\eta\Psi'_{j,k}(u+u\eta) dx=\int_{\R}u
\eta\Psi'_{j,k}dx+\int_{\R}u \eta^2\Psi'_{j,k}dx.
\end{eqnarray*}
Combining $J_1$-$J_4$ with (\ref{2.22})-(\ref{2.23}), we infer that
$\frac{d}{dt}I_j(t)$ can be written as
\begin{eqnarray}\label{2.24}
\frac{d}{dt}I_j(t)&=&-\frac{\dot{y}_j(t)}{2}\int_{\R}(u^2+u^2_x+\eta^2)\Psi'_{j,k}dx
+\frac{1}{2}\int_{\R}(uu^2_x+u\eta^2)
\Psi'_{j,k}dx\nonumber\\
&& +\int_{\R}u\Psi'_{j,k}(1-\partial^2_x)^{-1}(u^2+\frac{1}{2}u^2_x
+\frac{1}{2}\eta^2)dx+\int_{\R}u\Psi'_{j,k}(1-\partial^2_x)^{-1}\eta dx\nonumber\\
&\triangleq&-\frac{\dot{y}_j(t)}{2}\int_{\R}(u^2+u^2_x+\eta^2)\Psi'_{j,k}dx+J'_1+J'_2+J'_3.
\end{eqnarray}
To estimate $J'_1$, we divide $\R$ into two regions
$D_j\triangleq[x_{j-1}(t)+\frac{L}{4},x_j(t)-\frac{L}{4}]$ and its
complement $D^c_j$. For $x\in D_j$, from (\ref{2.2}), we obtain
\begin{eqnarray}\label{2.25}
\|u(t)\|_{L^\infty_{D_j}}
&\leq&\sum_{i=1}^N\|\varphi_{c_i}(\cdot-x_i(t))\|_{L^\infty_{D_j}}
+\|u-\sum_{i=1}^N\varphi_{c_i}(\cdot-x_i(t))\|_{L^\infty_{D_j}}\nonumber\\
&\leq& O(e^{-\sigma_0L})+O(\alpha).
\end{eqnarray}
On the other hand, for $x\in D^c_j$ and by (\ref{2.4}), we find
\begin{eqnarray}\label{2.26}
\big|x-y_j(t)\big|\geq\frac{x_j(t)-x_{j-1}(t)}{2}
-\frac{L}{4}\geq\frac{L}{8}+\sigma_0t.
\end{eqnarray}
Thus for $\alpha_0>0$ small enough and $L_0>0$ large enough, by
(\ref{2.25})-(\ref{2.26}), we obtain
\begin{eqnarray}\label{2.27}
J'_1&=&\frac{1}{2}\int_{\R}(uu^2_x+u\eta^2) \Psi'_{j,k}dx\nonumber\\
&=&\frac{1}{2}\int_{D^c_j}(uu^2_x+u\eta^2)
\Psi'_{j,k}dx+\frac{1}{2}\int_{D_j}(uu^2_x+u\eta^2) \Psi'_{j,k}dx\nonumber\\
&\leq&\frac{1}{2}\|u(t)\|_{L^\infty}\sup\limits_{x\in
D^c_j}\big|\Psi'_{j,K}(x-y_j(t))\big|\int_{\R}(u^2_x+\eta^2)
dx+\frac{1}{2}\|u(t)\|_{L^\infty_{D_j}}\int_{D_j}(u^2_x+\eta^2)
\Psi'_{j,k}dx\nonumber\\
&\leq&C\|\vec{u}^t_0\|^3_{X}e^{-\frac{(\sigma_0t+L/8)}{K}}+
\frac{\sigma_0}{8}\int_\R (u^2_x+\eta^2) \Psi'_{j,K}dx.
\end{eqnarray}
To bound $J'_2$, for $x\in D^c_j$, by the Young inequality, we get
\begin{eqnarray}\label{2.28}
&&\int_{D^c_j}u \Psi'_{j,K}
(1-\partial^2_x)^{-1}(u^2+\frac{1}{2}u^2_x
+\frac{1}{2}\eta^2) dx\nonumber\\
&\leq&\|u(t)\|_{L^\infty}\sup\limits_{x\in
D^c_j}\big|\Psi'_{j,K}(x-y_j(t))\big|\int_{\R}p\ast(u^2+u^2_x
+\eta^2)dx\nonumber\\
&=&\frac{1}{2}\|u(t)\|_{L^\infty}\sup\limits_{x\in
D^c_j}\big|\Psi'_{j,K}(x-y_j(t))\big|\int_{\R}e^{-|x|}\ast(u^2+u^2_x
+\eta^2)dx\nonumber\\
&\leq&C\|\vec{u}^t_0\|^3_{X}e^{-\frac{(\sigma_0t+L/8)}{K}},
\end{eqnarray}
where $p(x)\triangleq\frac{1}{2}e^{-|x|}$ is the Green function of
$(1-\partial_x^2)^{-1}.$ On the other hand, since by the definition
of $\Psi$, we find
\begin{eqnarray*}
(1-\partial^2_x)\Psi'_{j,K}=\Psi'_{j,K}-\frac{1}{K^3}\Psi'''(\frac{x-y_j(t)}{K})
\geq(1-\frac{10}{K^2})\Psi'_{j,K},
\end{eqnarray*}
hence, for $K\geq 5$,
\begin{eqnarray}\label{2.29}
(1-\partial^2_x)^{-1}\Psi'_{j,K}\leq(1-\frac{10}{K^2})^{-1}\Psi'_{j,K}.
\end{eqnarray}
For $x\in D_j$, noting that $\Psi'_{j,K}$ and $u^2+\frac{1}{2}u^2_x
+\frac{1}{2}\eta^2$ are non-negative, by (\ref{2.25}) and (\ref{2.29}), for
$\alpha_0$ small enough and $L_0$ large enough, we have
\begin{eqnarray}\label{2.30}
\int_{D_j}u \Psi'_{j,K} (1-\partial^2_x)^{-1}(u^2+\frac{1}{2}u^2_x
+\frac{1}{2}\eta^2) dx &\leq&\|u(t)\|_{L^\infty_{D_j}} \int_{D_j}
\Psi'_{j,K}(1-\partial^2_x)^{-1}(u^2+u^2_x
+\eta^2)dx\nonumber\\
&\leq&\|u\|_{L^\infty_{D_j}}\int_{\R} (u^2+u^2_x
+\eta^2)(1-\partial^2_x)^{-1}\Psi'_{j,K}dx\nonumber\\
&\leq& \frac{\sigma_0}{8}\int_\R (u^2+u^2_x+\eta^2) \Psi'_{j,K}dx.
\end{eqnarray}
It thus remains to estimate $J'_3$. For $s>0$ to be chosen later, by
the Cauchy-Schwarz inequality, we obtain
\begin{eqnarray*}
J'_3&=&\int_{\R}u\Psi'_{j,k}(1-\partial^2_x)^{-1}\eta
dx=\int_{\R}u\Psi'_{j,k} p\ast \eta dx \\
&\leq&
\frac{1}{2s}\int_{\R}u^2\Psi'_{j,k}dx+\frac{s}{2}\int_{\R}\Psi'_{j,k}
(p\ast \eta)^2 dx\leq
\frac{1}{2s}\int_{\R}u^2\Psi'_{j,k}dx+\frac{s}{2}J'_{31}.
\end{eqnarray*}
For $J'_{31}$, by H\"{o}lder inequality, we firstly get
\begin{eqnarray*}
(p\ast \eta)^2(x)&=&\big(\frac{1}{2}\int_{\R}e^{-|x-y|}\eta(y)dy\big)^2\\
&\leq& \frac{1}{4}\big(\int_{\R}e^{-|x-y|}dy\big)\cdot \big(\int_{\R}e^{-|x-y|}\eta^2(y)dy\big)
=\frac{1}{2}(p\ast \eta^2)(x).
\end{eqnarray*}
Thus, we compute the term $J'_{31}$ by the above inequality, the Minkowski inequality and (\ref{2.29}) as
\begin{eqnarray*}
J'_{31}&=&\int_{\R}\Psi'_{j,k}
(p\ast \eta)^2 dx \leq \frac{1}{2}\int_{\R}\Psi'_{j,k} (p\ast \eta^2) dx\\
&=&\frac{1}{2}\int_{\R}\Psi'_{j,k}(t,x)\big(\int_{\R}\frac{1}{2}e^{-|x-y|}\eta^2(y)dy\big)dx
\leq \frac{1}{2}\int_{\R}\eta^2(y)\int_{\R}\frac{1}{2}e^{-|y-x|}\Psi'_{j,k}(t,x)dxdy\\
&=&\frac{1}{2}\int_{\R}\eta^2(y)(1-\partial^2_x)^{-1}\Psi'_{j,K}dy\leq\frac{1}{2}(1-\frac{10}{K^2})^{-1}\int_{\R}\eta^2(y)\Psi'_{j,K}dy.
\end{eqnarray*}
Hence, we have
\begin{eqnarray}\label{2.31}
J'_3\leq
\frac{1}{2s}\int_{\R}u^2\Psi'_{j,k}dx+\frac{s}{4}(1-\frac{10}{K^2})^{-1}\int_{\R}\eta^2\Psi'_{j,K}dx.
\end{eqnarray}
For $\alpha_0>0$ small enough and $L_0>0$ large enough both
depending only on $\sigma_0>0$, it follows from (\ref{2.3}) that
\begin{eqnarray*}
-\frac{\dot{y}_j(t)}{2}=-\frac{\dot{x}_{j-1}(t)-c_{j-1}}{4}-\frac{\dot{x}_{j}(t)-c_{j}}{4}-\frac{c_{j-1}+c_j}{4}\leq-\frac{c_1+\sigma_0}{2}.
\end{eqnarray*}
Therefore, plugging (\ref{2.27})-(\ref{2.28}) and (\ref{2.30})-(\ref{2.31}) into (\ref{2.24}), and
by the above inequality, we derive that
\begin{eqnarray*}
\frac{d}{dt}I_j(t)&\leq&
C\|\vec{u}^t_0\|^3_{X}e^{-\frac{(\sigma_0t+L/8)}{K}}-\frac{\sigma_0}{4}\int_\R
(u^2+u^2_x+\eta^2)
\Psi'_{j,K}dx\nonumber\\
&&+(-\frac{c_1}{2}+\frac{1}{2s})\int_\R u^2
\Psi'_{j,K}dx+\big(-\frac{c_1}{2}+\frac{s}{4}(1-\frac{10}{K^2})^{-1}\big)\int_{\R}\eta^2\Psi'_{j,K}dx.
\end{eqnarray*}
Since $c_1>1$ and $K\geq5$, we then take $s>0$, such that
\begin{eqnarray*}
-\frac{c_1}{2}+\frac{1}{2s}\leq0\ \  \mbox{and}\
-\frac{c_1}{2}+\frac{s}{4}(1-\frac{10}{K^2})^{-1}\leq0.
\end{eqnarray*}
In this way we obtain
\begin{eqnarray*}
\frac{d}{dt}I_j(t)\leq
C\|\vec{u}^t_0\|^3_{X}e^{-\frac{(\sigma_0t+L/8)}{K}}-\frac{\sigma_0}{4}\int_\R
(u^2+u^2_x+\eta^2) \Psi'_{j,K}dx.
\end{eqnarray*}
Then the almost monotonicity property (\ref{2.19}) can be obtained by
integrating the above inequality from $0$ to $t$. This completes the
proof of Lemma \ref{Lem2}.
\end{proof}

\subsection{Local coercivity}

In this subsection, we present a local coerciveness inequality which
is crucial to our proof of the stability result. First, we recall
that the Hessian operator $H_c$ of $cE-F$ around a solitary wave
$\vec{\varphi}^t_c=(\varphi_c,\xi_c)^t$ is given by \cite{Z-L}
\begin{equation*}H_c= c
E''(\vec{\varphi}_c^t)-F''(\vec{\varphi}_c^t)= \left(
\begin{array}{cc}
 L_c& -(1+\xi_c)\\
-(1+\xi_c) &c-\varphi_c\\
\end{array}
\right),\end{equation*} where $L_c\triangleq
-\partial_x((c-\varphi_c)\partial_x)-3\varphi_c+\partial^2_x\varphi_{c}+c.$
Using $\xi_c=\frac{\varphi_c}{c-\varphi_c}$, we have
\begin{equation*}H_c= \left(
\begin{array}{cc}
 L_c& -\frac{c}{c-\varphi_c}\\
-\frac{c}{c-\varphi_c} &c-\varphi_c\\
\end{array}
\right).\end{equation*}
\begin{lemma2}\label{Lem3}
There exist $\delta,C_\delta,C>0$ depending only on $c_1>1$, such
that for all $c\geq c_1$, $\Theta(x)\in C^2(\R)>0$ and
$\vec{\psi}^t=(\psi,\omega)^t\in X,$ satisfying
\begin{eqnarray}\label{2.32}
\big|\big<\sqrt{\Theta}\vec{\psi}^t,\big((1-\partial^2_x)\varphi_c,\xi_c
\big)^t\big>_{L^2\times L^2}\big|
+\big|\big<\sqrt{\Theta}\vec{\psi}^t,\big((1-\partial^2_x)\partial_x\varphi_c,\partial_x\xi_c
\big)^t \big>_{L^2\times L^2}\big|\leq \delta\|\vec{\psi}^t\|_X,
\end{eqnarray}
and
\begin{eqnarray}\label{2.33}
\big|\frac{(\Theta')^2}{4\Theta}\big|+c\big|\Theta'\big|+\big|\frac{\Theta''}{2}\big|
\leq\min\{\frac{1}{4},\frac{C_\delta}{4c}\}\Theta.
\end{eqnarray}
Then, we have
\begin{eqnarray}\label{2.34}
\Lambda&\triangleq&\int_{\R}\Theta\big((c-\varphi_c)(\partial_x\psi)^2+(-3\varphi_c+\partial^2_x\varphi_c+c)\psi^2\big)
+\partial_x\varphi_c\Theta'\psi^2-2\Theta\frac{c}{c-\varphi_c}\psi\omega\nonumber\\
&&+\Theta(c-\varphi_c)\omega^2dx \geq C
\int_{\R}\Theta(\psi^2+(\partial_x\psi)^2+\omega^2)dx.
\end{eqnarray}
\end{lemma2}
\begin{proof}
We directly calculate that
\begin{eqnarray}\label{2.35}
&&\big< H_c\sqrt{\Theta}\vec{\psi}^t,
\sqrt{\Theta}\vec{\psi}^t\big>_{L^2\times L^2}= \Big\langle
H_c\sqrt{\Theta}\left(
\begin{array}{c}
\psi\\
\omega\\
\end{array}
\right), \sqrt{\Theta}\left(
\begin{array}{c}
\psi\\
\omega\\
\end{array}
\right)\Big\rangle_{L^2\times L^2} \nonumber\\
&=&\Big\langle \left(
\begin{array}{c}
\big(-\partial_x((c-\varphi_c)\partial_x)-3\varphi_c+\partial^2_x\varphi_{c}+c\big)\sqrt{\Theta}\psi
-\frac{c}{c-\varphi_c}\sqrt{\Theta}\omega\\
-\frac{c}{c-\varphi_c}\sqrt{\Theta}\psi+(c-\varphi_c)\sqrt{\Theta}\omega\\
\end{array}
\right),\left(
\begin{array}{c}
 \sqrt{\Theta}\psi\\
 \sqrt{\Theta}\omega\\
\end{array}
\right)\Big\rangle_{L^2\times L^2}\nonumber\\
&=&\int_\R
\Theta\big((c-\varphi_c)(\partial_x\psi)^2+(-3\varphi_c+\partial^2_x\varphi_{c}+c)\psi^2\big)dx
-2\int_\R\Theta\frac{c}{c-\varphi_c}\psi\omega dx\nonumber\\
&&+\int_\R\Theta(c-\varphi_c)\omega^2dx+\int_\R\frac{c-\varphi_c}{4}\frac{(\Theta')^2}{\Theta}\psi^2dx+\int_\R(c-\varphi_c)\Theta'\psi\partial_x\psi
dx\nonumber\\
&=&\int_{\R}\Theta\big((c-\varphi_c)(\partial_x\psi)^2+(-3\varphi_c+\partial^2_x\varphi_c+c)\psi^2\big)
+\partial_x\varphi_c\Theta'\psi^2-2\Theta\frac{c}{c-\varphi_c}\psi\omega\nonumber\\
&&+\Theta(c-\varphi_c)\omega^2dx+\int_\R\big((c-\varphi_c)(\frac{(\Theta')^2}{4\Theta}-\frac{\Theta''}{2})-\frac{1}{2}\partial_x\varphi_c\Theta'\big)\psi^2dx\nonumber\\
&=&\Lambda+\int_\R\big((c-\varphi_c)(\frac{(\Theta')^2}{4\Theta}-\frac{\Theta''}{2})-\frac{1}{2}\partial_x\varphi_c\Theta'\big)\psi^2dx,
\end{eqnarray}
and
\begin{eqnarray}\label{2.36}
\|\sqrt{\Theta}\vec{\psi}^t\|_{X}^2&=&\int_\R(\sqrt{\Theta}\psi)^2+(\partial_x(\sqrt{\Theta}\psi))^2+(\sqrt{\Theta}\omega)^2dx\nonumber\\
&=&\int_\R\Theta(\psi^2+(\partial_x\psi)^2+\omega^2)dx+\int_\R(\frac{(\Theta')^2}{4\Theta}-\frac{\Theta''}{2})\psi^2dx.
\end{eqnarray}
By the analysis on the spectrum of $H_c$ given in \cite{Z-L}, we can
easily deduce that there exist $\delta>0$ and $C_\delta>0$, such
that if for $c\geq c_1,$
\begin{eqnarray*}
\big|\big<\vec{\psi}^t,\big((1-\partial^2_x)\varphi_c,\xi_c
\big)^t\big>_{L^2\times L^2}\big|
+\big|\big<\vec{\psi}^t,\big((1-\partial^2_x)\partial_x\varphi_c,\partial_x\xi_c
\big)^t \big>_{L^2\times L^2}\big|\leq \delta\|\vec{\psi}^t\|_X,
\end{eqnarray*}
then
\begin{eqnarray*}
\big<H_c\vec{\psi}^t,\vec{\psi}^t\big>_{L^2\times L^2}\geq
C_\delta\|\vec{\psi}^t\|^2_X.
\end{eqnarray*}
Therefore, under the hypotheses (\ref{2.32})-(\ref{2.33}), we can derive (\ref{2.34})
from (\ref{2.35})-(\ref{2.36}) that
\begin{eqnarray*}
\Lambda&+&c\cdot \frac{C_\delta}{4c}\int_\R \Theta\psi^2dx\geq
\big<H_c\sqrt{\Theta}\vec{\psi}^t,
\sqrt{\Theta}\vec{\psi}^t\big>\geq C_\delta
\|\sqrt{\Theta}\vec{\psi}^t\|_{X}^2\\ &&\geq
C_\delta\int_\R\Theta(\psi^2+(\partial_x\psi)^2+\omega^2)dx
-\min\{\frac{1}{4},\frac{C_\delta}{4c}\}\int_\R \Theta\psi^2dx,
\end{eqnarray*}
where we used the fact that $\varphi_c\in [0,c-1]$ and
$\partial_x\varphi_c\in[-c+1,c-1]$. This completes the proof of
Lemma \ref{Lem3}.
\end{proof}

\section{Proof of the orbital stability}\label{Sec3}
\newtheorem {remark3}{Remark}[section]
\newtheorem{theorem3}{Theorem}[section]
\newtheorem{lemma3}{Lemma}[section]
Based on the series of lemmas in Section \ref{Sec2}, we will complete the proof
of the orbital stability of the train of $N$ solitary waves for the
system (\ref{2-CH}). By the continuity of $\vec{u}^t(t)$ in $H^s\times
H^{s-1}\hookrightarrow X$, with $s\geq 2$, to prove Theorem \ref{Th1}, it
is sufficient to show that there exist $A_0,\ \gamma_0,\ L_0,\
\varepsilon_0>0$ such that for all $ L>L_0$ and
$0<\varepsilon<\varepsilon_0$, if $\vec{u}^t_0$ satisfies (\ref{1.6}) and
for some $0<t_0<T$, with $0<T\leq +\infty,$
\begin{eqnarray}\label{3.1}
\vec{u}^t(t)\in U(A_0(\sqrt\varepsilon+e^{-\gamma_0L}),\frac{L}{2}),
\quad \mbox{for\ all}\  \ t\in [0,t_0],
\end{eqnarray}
then
\begin{eqnarray}\label{3.2}
\vec{u}^t(t_0)\in
U(\frac{A_0}{2}(\sqrt\varepsilon+e^{-\gamma_0L}),\frac{2L}{3}).
\end{eqnarray}
Therefore, we will conclude the proof of Theorem \ref{Th1}, if we prove
the result (\ref{3.2}) under the assumption (\ref{3.1}) for some $ L>L_0$ and
$0<\varepsilon<\varepsilon_0$ with $A_0,\ \gamma_0,\ L_0,\
\varepsilon_0>0$ to be specified later.

For simplicity, we set $\vec{u}^t=\vec{u}^t(t_0),$
$X=(x_1,...,x_N)=(x_1(t_0),...,x_N(t_0)).$ For $i=1,...,N$, we
define $\vec{\psi}_i^t=(\psi_i,\omega_i)^t\in X$ by
\begin{eqnarray}\label{3.3}
\vec{u}^t=(1+a_i)\vec{R}^t_X+\vec{\psi}_i^t, \quad\quad
\big<E_i'(\vec{R}^t_X), \vec{\psi}_i^t\big>_{L^2\times L^2}=0,
\end{eqnarray}
where
$\vec{R}^t_X=(R_X,S_X)^t=\sum\limits_{i=1}\limits^N\vec{R}^t_i=
\big(\sum\limits_{i=1}\limits^NR_i(\cdot),\sum\limits_{i=1}\limits^NS_i(\cdot)\big)^t
=\big(\sum\limits_{i=1}\limits^N\varphi_{c_i}(\cdot-x_i),\sum\limits_{i=1}\limits^N\xi_{c_i}(\cdot-x_i)\big)^t.$
Since, by (\ref{2.17})-(\ref{2.18}) and (\ref{1.2}), we have
\begin{eqnarray}\label{3.4}
&&\big<(E^t_i)'(\vec{R}^t_{X}),\vec{R}^t_{X}\big>_{L^2\times
L^2}=\big<E'(\vec{\varphi}^t_{c_i}),\vec{\varphi}^t_{c_i}\big>_{L^2\times
L^2}+O(e^{-\sigma_0L})\nonumber\\&=&
\big\|\vec{\varphi}^t_{c_i}\big\|^2_{X}+O(e^{-\sigma_0L})>\frac{1}{2}\big\|\vec{\varphi}^t_{c_i}\big\|^2_{X}.
\end{eqnarray}
Thus, the functions $\vec{\psi}_i^t$ are well defined. Then, we set
$\vec{v}^t=\big(v(t,x),\zeta(t,x)\big)^t=\vec{u}^t-\vec{R}^t_X$ and
in the sequel suppose that
\begin{eqnarray}\label{3.5}
\big\|\vec{v}^t\big\|_{X}\geq
\sqrt{\varepsilon}+e^{-\frac{\sigma_0L}{2}},
\end{eqnarray}
otherwise we complete our proof with $A_0=2$ and
$\gamma_0=\frac{\sigma_0}{2}$. We break our proof into three steps.
\begin{proof}
$Step\ 1:$ In the first step, for all $i=1,...,N,$ we claim that the
following estimates on $a_i:$
\begin{eqnarray}\label{3.6}
\big|a_i\big|\leq O(\big\|\vec{v}^t\big\|^2_{X}),
\end{eqnarray} hold true.

Indeed, according to the definitions (\ref{2.16}) of $I_i,$ $E^t_i$ and
$F^t_i,$ we have
\begin{eqnarray}\label{3.7}
I_i(t,\vec{u}^t)=\sum\limits_{j=i}\limits^NE^t_j(\vec{u}^t),\ \
\mbox{for} \ \ i=2,...,N,  \ \
E(\vec{u}^t)=\sum\limits_{j=1}\limits^NE^t_j(\vec{u}^t)\ \
\mbox{and}\ \
F(\vec{u}^t)=\sum\limits_{j=1}\limits^NF^t_j(\vec{u}^t) .
\end{eqnarray}
Using the exponentially asymptotic behavior of $\vec{\varphi}_c^t$,
and by (\ref{2.17})-(\ref{2.18}), one can easily find that
\begin{eqnarray}\label{3.8}
E^t_j\big((R_j,S_j)^t\big)=E(\vec{\varphi}^t_{c_j})+O(e^{-\sigma_0L})\
\ \mbox{and}\ \ E^t_j\big((R_k,S_k)^t\big)\leq O(e^{-\sigma_0L})\ \
\mbox{for}\ \ j\neq k.
\end{eqnarray}
Hence, by Taylor formula, (\ref{3.3}), (\ref{3.5}) and (\ref{3.8}), we obtain
\begin{eqnarray}\label{3.9}
\sum\limits_{j=1}\limits^NE^t_j(\vec{u}^t)&=&\sum\limits_{j=1}\limits^NE^t_j(\vec{R}^t_{X})+\sum\limits_{j=1}\limits^N\big<(E^t_j)'(\vec{R}^t_{X}),\vec{v}^t\big>_{L^2\times
L^2}+O(\|\vec{v}^t\|^2_{X})\nonumber\\
&=&\sum\limits_{j=1}\limits^NE(\vec{\varphi}^t_{c_j})+\sum\limits_{j=1}\limits^Na_j\big<(E^t_j)'(\vec{R}^t_{X}),\vec{R}^t_{X}\big>_{L^2\times
L^2}+O(\|\vec{v}^t\|^2_{X}).
\end{eqnarray}
Since $\vec{u}_0^t$ satisfies (\ref{1.6}), on account of the conservation
laws (\ref{1.5}) of $E$ and $F$, we get
\begin{eqnarray}\label{3.10}
E(\vec{u}^t)=E(\vec{u}_0^t)=\sum\limits_{j=1}\limits^NE(\vec{\varphi}^t_{c_j})+O(e^{-\sigma_0L})+O(\varepsilon),
\end{eqnarray}
and
\begin{eqnarray}\label{3.11}
F(\vec{u}^t)=F(\vec{u}_0^t)=\sum\limits_{j=1}\limits^NF(\vec{\varphi}^t_{c_j})+O(e^{-\sigma_0L})+O(\varepsilon).
\end{eqnarray}
Thus, for all $i=1,...,N$, we deduce from (\ref{3.7}), (\ref{3.9})-(\ref{3.10}) and
(\ref{3.5}) that
\begin{eqnarray}\label{3.12}
\sum\limits_{j=i}\limits^Na_j\big<(E^t_j)'(\vec{R}^t_{X}),\vec{R}^t_{X}\big>_{L^2\times
L^2}\leq O(\|\vec{v}^t\|^2_{X}).
\end{eqnarray}
On the other hand, in a similar way, using Taylor formula, and by
(\ref{2.17})-(\ref{2.18}), we have
\begin{eqnarray*}\label{3.13}
F(\vec{u}^t)=\sum\limits_{i=1}\limits^NF^t_i(\vec{u}^t)=\sum\limits_{i=1}\limits^NF^t_i(\vec{R}^t_{X})+\sum\limits_{i=1}\limits^N\big<(F^t_i)'(\vec{R}^t_{X}),\vec{v}^t\big>_{L^2\times
L^2}+O(\|\vec{v}^t\|^2_{X}),
\end{eqnarray*}
and
\begin{eqnarray}\label{3.14}
F(\vec{u}^t)=\sum\limits_{i=1}\limits^NF(\vec{\varphi}^t_{c_i})+\sum\limits_{i=1}\limits^N\big<(F^t_i)'(\vec{R}^t_{X}),\vec{v}^t\big>_{L^2\times
L^2}+O(\|\vec{v}^t\|^2_{X})+O(e^{-\sigma_0L}).
\end{eqnarray}
Hence, by (\ref{3.11}), (\ref{3.14}), (\ref{3.5}) and (\ref{3.3}), we get
\begin{eqnarray}\label{3.15}
&&O(\|\vec{v}^t\|^2_{X})\nonumber\\&=&
\sum\limits_{i=1}\limits^N\big<(F^t_i)'(\vec{R}^t_{X}),\vec{v}^t\big>_{L^2\times
L^2}\nonumber\\
&=&\sum\limits_{i=1}\limits^N\big<(F^t_i)'(\vec{R}^t_{X})-c_i(E^t_i)'(\vec{R}^t_{X}),\vec{v}^t\big>_{L^2\times
L^2}+\sum\limits_{i=1}\limits^N\big<c_i(E^t_i)'(\vec{R}^t_{X}),\vec{v}^t\big>_{L^2\times
L^2}\nonumber\\
&=&\sum\limits_{i=1}\limits^N\big<(F^t_i)'(\vec{R}^t_{X})-c_i(E^t_i)'(\vec{R}^t_{X}),\vec{v}^t\big>_{L^2\times
L^2}+\sum\limits_{i=1}\limits^Nc_ia_i\big<(E^t_i)'(\vec{R}^t_{X}),\vec{R}^t_{X}\big>_{L^2\times
L^2}.
\end{eqnarray}
Since the solitary waves satisfy the identity
$cE'(\vec{\varphi}_c^t)-F'(\vec{\varphi}_c^t)=0$, by (\ref{2.17})-(\ref{2.18}),
then we obtain
\begin{eqnarray}\label{3.16}
\big\|(F^t_i)'(\vec{R}^t_{X})-c_i(E^t_i)'(\vec{R}^t_{X})\big\|_{X^\ast}\leq\big\|F'(\vec{\varphi}^t_{c_i})-c_iE'(\vec{\varphi}^t_{c_i})\big\|_{X^\ast}+O(e^{-\sigma_0L})\leq
O(e^{-\sigma_0L}).
\end{eqnarray}
Thus, by (\ref{3.15})-(\ref{3.16}), we infer that
\begin{eqnarray*}\label{3.17}
\sum\limits_{i=1}\limits^Nc_ia_i\big<(E^t_i)'(\vec{R}^t_{X}),\vec{R}^t_{X}\big>_{L^2\times
L^2}=O(\|\vec{v}^t\|^2_{X}).
\end{eqnarray*}
Then, using Abel transformation, we find that
\begin{eqnarray}\label{3.18}
&&\sum\limits_{i=1}\limits^N(c_i-c_{i-1})\sum\limits_{j=i}\limits^Na_j\big<(E^t_j)'(\vec{R}^t_{X}),\vec{R}^t_{X}\big>_{L^2\times
L^2}
+c_1\sum\limits_{j=1}\limits^Na_j\big<(E^t_j)'(\vec{R}^t_{X}),\vec{R}^t_{X}\big>_{L^2\times
L^2}\nonumber\\ &=&O(\|\vec{v}^t\|^2_{X}).
\end{eqnarray}
Combining (\ref{3.12}) with (\ref{3.18}), for all $i=1,...,N,$ it follows that
\begin{eqnarray*}\label{3.19}
\big|a_i\big<(E^t_i)'(\vec{R}^t_{X}),\vec{R}^t_{X}\big>_{L^2\times
L^2}\big|\leq O(\|\vec{v}^t\|^2_{X}).
\end{eqnarray*}
Therefore, by (\ref{3.4}), we prove our claim (\ref{3.6}). Moreover, for all
$i=1,...,N,$ we derive from our claim (\ref{3.6}) and (\ref{3.3}) that
\begin{eqnarray}\label{3.20}
\|\vec{v}^t\|_{X} \thicksim  \|\vec{\psi}_i^t\|_{X}.
\end{eqnarray}

$Step\ 2:$ In the second step, we will apply the local coerciveness
inequality (\ref{2.34}) in Lemma \ref{Lem3} to prove
\begin{eqnarray}\label{3.21}
\big<H_i^t(\vec{R}^t_{i})\vec{\psi}^t_{i},\vec{\psi}^t_{i}\big>_{L^2\times
L^2}\geq C E^t_i(\vec{\psi}^t_{i}),
\end{eqnarray}
where $\vec{R}^t_i(\cdot)=\big(R_i(\cdot),S_i(\cdot)\big)^t
=\big(\varphi_{c_i}(\cdot-x_i),\xi_{c_i}(\cdot-x_i)\big)^t$ and
$\vec{\psi}_i^t=(\psi_i,\omega_i)^t\in X$. And the operator $H_i^t$
is given in the following form (\ref{3.22}).

From the definitions of $
E^t_i\ \mbox{and} \ F^t_i$, we can explicitly compute the
variational derivatives as in \cite{Liu1}
\begin{equation*}
\left\{\begin{array}{ll}(E^t_i)'_\varphi=\varphi\Phi_i-\partial^2_x\varphi\Phi_i-\partial_x\varphi\partial_x\Phi_i,\\
(E^t_i)'_\xi=\xi \Phi_i,
\end{array}\right.
\end{equation*}
and
\begin{equation*}
\left\{\begin{array}{ll}(F^t_i)'_\varphi=\frac{3}{2}\varphi^2\Phi_i-\frac{1}{2}(\partial_x\varphi)^2\Phi_i
-\varphi\partial^2_x\varphi\Phi_i+\xi\Phi_i+\frac{1}{2}\xi^2\Phi_i-\varphi\partial_x\varphi\partial_x\Phi_i,\\
(F^t_i)'_\xi=\varphi \Phi_i+\varphi \xi \Phi_i.
\end{array}\right.
\end{equation*}
Hence, we have
\begin{equation*}
\left\{\begin{array}{lll}(E^t_i)''_{\varphi\varphi}=\Phi_i(1-\partial^2_x)-\partial_x\Phi_i\partial_x,\\
(E^t_i)''_{\varphi\xi}=(E^t_i)''_{\varphi\xi}=0,\\
(E^t_i)''_{\xi\xi}=\Phi_i,
\end{array}\right.
\end{equation*}
and
\begin{equation*}
\left\{\begin{array}{lll}(F^t_i)''_{\varphi\varphi}=(3\varphi-\partial_x\varphi\partial_x-\varphi\partial_x^2
-\partial_x^2\varphi)\Phi_i+(-\varphi\partial_x\Phi_i\partial_x-\partial_x\varphi\partial_x\Phi_i),\\
(F^t_i)''_{\varphi\xi}=(F^t_i)''_{\varphi\xi}=\Phi_i+\Phi_i\xi,\\
(F^t_i)''_{\xi\xi}=\varphi\Phi_i,
\end{array}\right.
\end{equation*}
Therefore, the linearized operator $H^t_i$ of $c(E^t_i)'-(F^t_i)'$
(or Hessian of $cE^t_i-F^t_i$) at $\vec{\varphi}_c^t$ can be given
as
\begin{eqnarray*}H^t_i(\vec{\varphi}_c^t)&=& c
(E^t_i)''(\vec{\varphi}_c^t)-(F^t_i)''(\vec{\varphi}_c^t)\\
&=& \left(
\begin{array}{cc}
 L_i^t& -\Phi_i(1+\xi_c)\\
-\Phi_i(1+\xi_c) &\Phi_i(c-\varphi_c)\\
\end{array}
\right),\end{eqnarray*} where $L_i^t\triangleq
-\partial_x(\Phi_i(c-\varphi_c)\partial_x)+\Phi_i(-3\varphi_c+\partial^2_x\varphi_{c}+c)+\partial_x\varphi\partial_x\Phi_i.$
By $\xi_c=\frac{\varphi_c}{c-\varphi_c}$, we have
\begin{equation}\label{3.22}
H^t_i(\vec{\varphi}_c^t)= \left(
\begin{array}{cc}
 L_i^t& -\frac{c\Phi_i}{c-\varphi_c}\\
-\frac{c\Phi_i}{c-\varphi_c} &\Phi_i(c-\varphi_c)\\
\end{array}
\right).
\end{equation}
According to Lemma \ref{Lem3}, to prove (\ref{3.21}), we only need to verify that
$\vec{\psi}_i^t$ satisfies the condition (\ref{2.32}) with $\Theta=\Phi_i$
and $\vec{\varphi}^t_c=\vec{R}^t_i$. Indeed, using the fact
$\vec{\psi}^t_{i}=(\psi_{i},\omega_{i})^t=\vec{v}^t-a_i\vec{R}^t_X=(v-a_iR_X,\zeta-a_iS_X)^t$,
by (\ref{2.1}), we get
\begin{eqnarray*}
&&\big<\sqrt{\Phi_i}\vec{\psi}^t_{i},\big((1-\partial^2_x)\partial_xR_{i},\partial_xS_{i}\big)^t\big>_{L^2\times L^2}\nonumber\\
&=&\big<\vec{\psi}^t_{i},\big((1-\partial^2_x)\partial_xR_{i},\partial_xS_{i}\big)^t\big>_{L^2\times L^2}
+\big<(\sqrt{\Phi_i}-1)\vec{\psi}^t_{i},\big((1-\partial^2_x)\partial_xR_{i},\partial_xS_{i}\big)^t\big>_{L^2\times L^2}\nonumber\\
&=&\big<\psi_{i},(1-\partial^2_x)\partial_xR_{i}\big>_{L^2}+\big<\omega_{i},\partial_xS_{i}\big>_{L^2}+\big<(\sqrt{\Phi_i}-1)\psi_{i},(1-\partial^2_x)\partial_xR_{i}\big>_{L^2}\nonumber\\
&&+\big<(\sqrt{\Phi_i}-1)\omega_{i},\partial_xS_{i}\big>_{L^2}\nonumber\\
&=&-a_i\big<R_X,(1-\partial^2_x)\partial_xR_{i}\big>_{L^2}-a_i\big<S_X,\partial_xS_{i}\big>_{L^2}+\big<(\sqrt{\Phi_i}-1)\psi_{i},(1-\partial^2_x)\partial_xR_{i}\big>_{L^2}\\
&&+\big<(\sqrt{\Phi_i}-1)\omega_{i},\partial_xS_{i}\big>_{L^2}.
\end{eqnarray*}
Thus, gathering (\ref{2.17}), (\ref{3.1}), (\ref{3.6}) and (\ref{3.20}), we deduce that
\begin{eqnarray*}
\big|\big<\sqrt{\Phi_i}\vec{\psi}^t_{i},\big((1-\partial^2_x)\partial_xR_{i},\partial_xS_{i}\big)^t\big>_{L^2\times
L^2}\big|\leq \big(
O(\sqrt{\varepsilon})+O(e^{-\sigma_0L})+O(e^{-\gamma_0L})\big)\|\vec{\psi}^t_{i}\|_{X}.
\end{eqnarray*}
In a similar manner as above and by (\ref{3.3}), we then obtain
\begin{eqnarray*}
&&\big<\sqrt{\Phi_i}\vec{\psi}^t_{i},\big((1-\partial^2_x)R_{i},S_{i}\big)^t\big>_{L^2\times L^2}\\
&=&\big<\vec{\psi}^t_{i},\big((1-\partial^2_x)R_{i},S_{i}\big)^t\big>_{L^2\times L^2}+\big<(\sqrt{\Phi_i}-1)\vec{\psi}^t_{i},\big((1-\partial^2_x)R_{i},S_{i}\big)^t\big>_{L^2\times L^2}\nonumber\\
&=&\big<\psi_{i},(1-\partial^2_x)R_{i}\big>_{L^2}+\big<\omega_{i},S_{i}\big>_{L^2}+\big<(\sqrt{\Phi_i}-1)\psi_{i},(1-\partial^2_x)R_{i}\big>_{L^2}+\big<(\sqrt{\Phi_i}-1)\omega_{i},S_{i}\big>_{L^2}\nonumber\\
&=&\big<(E^t_i)'(\vec{R}^t_{X}),\vec{\psi}^t_{i}\big>_{L^2\times L^2}+\int_\R(1-\Phi_i)\big((\psi_iR_i+\partial_x\psi_i\partial_xR_i)+\omega_iS_i\big)dx\\
&&-\sum\limits_{j=1,j\neq i}\limits^N\int_\R
\Phi_i\big((\psi_iR_j+\partial_x\psi_i\partial_xR_j)+\omega_iS_j\big)dx+\big<(\sqrt{\Phi_i}-1)\psi_{i},(1-\partial^2_x)R_{i}\big>_{L^2}\\
&&+\big<(\sqrt{\Phi_i}-1)\omega_{i},S_{i}\big>_{L^2} \\
&\leq&\big(
O(\sqrt{\varepsilon})+O(e^{-\sigma_0L})+O(e^{-\gamma_0L})\big)\|\vec{\psi}^t_{i}\|_{X}.
\end{eqnarray*}

$Step\ 3:$ In the last step, we prove that there exists $C>0$
independent of $A_0$ and $\gamma_0$ such that for $L\geq L_0$, with
$L_0$ large enough,
\begin{eqnarray*}
\big\|\vec{v}^t\big\|^2_{X}=
\big\|\vec{u}^t-\sum_{i=1}^N\vec{\varphi}^t_{c_i}(\cdot-x_i)\big\|^2_{X}
\leq C(\varepsilon+e^{-\sigma_0L}).
\end{eqnarray*}
Then, we can conclude the proof of Theorem \ref{Th1} with
$\gamma_0=\frac{\sigma_0}{2}$ and $A_0=2\sqrt{C}.$

By the monotonicity estimates (\ref{2.19}) and (\ref{1.6}), we have
\begin{eqnarray}\label{3.23}
I_i(t_0,\vec{u}^t)\leq I_i(0,\vec{u}_0^t)+O(e^{-\sigma_0L})&\leq&
\sum\limits_{j=i}\limits^NE_j^0\big(\vec{\varphi}^t_{c_j}(\cdot-x_j(0))\big)+O(e^{-\sigma_0L})+O(\varepsilon)\nonumber\\
&\leq&
\sum\limits_{j=i}\limits^NE(\vec{\varphi}^t_{c_j})+O(e^{-\sigma_0L})+O(\varepsilon).
\end{eqnarray}
Using Abel transformation, in view of the conservation law $E$, then
by (\ref{3.7}), (\ref{3.23}) and (\ref{3.9}), we deduce that
\begin{eqnarray}\label{3.24}
\sum\limits_{i=1}\limits^Nc_iE^t_i(\vec{u}^t)&=&\sum\limits_{i=2}\limits^N(c_i-c_{i-1})\sum\limits_{j=i}\limits^NE^t_j(\vec{u}^t)+c_1
\sum\limits_{j=1}\limits^NE^t_j(\vec{u}^t)\nonumber\\
&=&\sum\limits_{i=2}\limits^N(c_i-c_{i-1})I_i(t_0,\vec{u}^t(t_0))+c_1E(\vec{u}^t(t_0))\nonumber\\
&\leq& \sum\limits_{i=2}\limits^N(c_i-c_{i-1})I_i(0,\vec{u}^t_0)+c_1E(\vec{u}^t_0)+O(e^{-\sigma_0 L})\nonumber\\
&\leq&\sum\limits_{i=2}\limits^N(c_i-c_{i-1})\sum\limits_{j=i}\limits^NE(\vec{\varphi}^t_{c_j})+c_1
\sum\limits_{j=1}\limits^NE(\vec{\varphi}^t_{c_j})+O(e^{-\sigma_0 L})+O(\varepsilon)\nonumber\\
&\leq&
\sum\limits_{i=1}\limits^Nc_iE(\vec{\varphi}^t_{c_i})+O(e^{-\sigma_0
L})+O(\varepsilon)\nonumber\\
&\leq&\sum\limits_{i=1}\limits^Nc_iE^t_i(\vec{R}^t_{X})+O(e^{-\sigma_0
L})+O(\varepsilon).
\end{eqnarray}
On the other hand, using Taylor formula, we derive from (\ref{3.3}) and
(\ref{3.16}) that
\begin{eqnarray}\label{3.25}
&&\sum\limits_{i=1}\limits^N\big(c_iE^t_i(\vec{u}^t)-F^t_i(\vec{u}^t)\big)\nonumber\\&=&\sum\limits_{i=1}\limits^N\big(c_iE^t_i(\vec{R}^t_{X})-F^t_i(\vec{R}^t_{X})\big)
+\frac{1}{2}\sum\limits_{i=1}\limits^N\big<H^t_i(\vec{R}^t_{X})\vec{\psi}_i^t,\vec{\psi}_i^t\big>_{L^2\times
L^2}+\sum\limits_{i=1}\limits^Na_i\big<H^t_i(\vec{R}^t_{X})\vec{R}^t_{X},\vec{\psi}_i^t\big>_{L^2\times
L^2}\nonumber\\
&&+\sum\limits_{i=1}\limits^N\frac{a_i^2}{2}\big<H^t_i(\vec{R}^t_{X})\vec{R}^t_{X},\vec{R}^t_{X}\big>_{L^2\times
L^2}+o(\|\vec{v}^t\|^2_{X})+O(e^{-\sigma_0 L}).
\end{eqnarray}
Combining (\ref{3.24}) with (\ref{3.25}), we obtain
\begin{eqnarray}\label{3.26}
&&F(\vec{R}^t_{X})-F(\vec{u}^t)\nonumber\\
&\geq&\frac{1}{2}\sum\limits_{i=1}\limits^N\big<H^t_i(\vec{R}^t_{X})\vec{\psi}_i^t,\vec{\psi}_i^t\big>_{L^2\times
L^2}+\sum\limits_{i=1}\limits^Na_i\big<H^t_i(\vec{R}^t_{X})\vec{R}^t_{X},\vec{\psi}_i^t\big>_{L^2\times
L^2}\nonumber\\
&&+\sum\limits_{i=1}\limits^N\frac{a_i^2}{2}\big<H^t_i(\vec{R}^t_{X})\vec{R}^t_{X},\vec{R}^t_{X}\big>_{L^2\times
L^2}+o(\|\vec{v}^t\|^2_{X})+O(e^{-\sigma_0 L})+O(\varepsilon)
\end{eqnarray}
Since, from (\ref{3.20}), (\ref{3.22}) and (\ref{3.6}), we get
\begin{eqnarray*}
\sum\limits_{i=1}\limits^Na_i\big<H^t_i(\vec{R}^t_{X})\vec{R}^t_{X},\vec{\psi}_i^t\big>_{L^2\times
L^2}=O(\|\vec{v}^t\|^3_{X})+O(e^{-\sigma_0 L})+O(\varepsilon),
\end{eqnarray*}
and we can easily infer that
\begin{eqnarray*}
\big<H^t_i(\vec{R}^t_{X})\vec{\psi}_i^t,\vec{\psi}_i^t\big>_{L^2\times
L^2}=\big<H_i^t(\vec{R}^t_{i})\vec{\psi}^t_{i},\vec{\psi}^t_{i}\big>_{L^2\times
L^2}+O(e^{-\sigma_0 L})+O(\|\vec{v}^t\|^3_{X}),
\end{eqnarray*}
Hence, by (\ref{3.21}), (\ref{3.26}) and the above two equalities,
\begin{eqnarray}\label{3.27}
\sum\limits_{i=1}\limits^NE^t_i(\vec{\psi}^t_{i})\leq
o(\|\vec{v}^t\|^2_{X})+O(e^{-\sigma_0 L})+O(\varepsilon).
\end{eqnarray}
Then, by (\ref{3.3}) and (\ref{3.6}), we get
\begin{eqnarray}\label{3.28}
\sum\limits_{i=1}\limits^NE^t_i(\vec{\psi}^t_{i})=\sum\limits_{i=1}\limits^NE^t_i(\vec{v}^t)+O(\|\vec{v}^t\|^3_{X})
=E(\vec{v}^t)+O(\|\vec{v}^t\|^3_{X})\geq \frac{1}{2}E(\vec{v}^t).
\end{eqnarray}
Therefore, gathering (\ref{3.26})-(\ref{3.28}) and (\ref{3.11}), we complete the proof
of Theorem \ref{Th1}.
\end{proof}

\bigskip
\noindent\textbf{Acknowledgments} The author thanks the anonymous
referee for helpful suggestions and comments. The work is supported by the National
Natural Science Foundation of China (No.11426212).

\end{document}